\documentclass[runningheads]{llncs}
\usepackage[T1]{fontenc}
\usepackage[utf8]{inputenc}
\usepackage{subcaption}
\usepackage{amsmath}
\usepackage{tikz}
\usepackage{multirow}
\usepackage{multicol}
\usepackage{booktabs}
\usepackage{pgfplots}
\usepackage{hhline}
\usepackage{tabularx}
\usepackage{adjustbox}
\usepackage{float}
\usepackage{amssymb}
\usepackage{listings}

\usetikzlibrary{positioning,arrows,shapes,shapes.geometric}
\usepackage{graphicx}

\begin{document}

\title{Observe Locally, Classify Globally: Using GNNs to Identify Sparse Matrix Structure}
\titlerunning{Using GNNs to Identify Sparse Matrix Structure}

\author{Khaled Abdelaal\inst{1}\orcidID{0000-0001-8963-7523} \and
Richard Veras\inst{1}\orcidID{0000-0003-2633-3391}}
\authorrunning{K. Abdelaal and R. Veras}

\institute{University of Oklahoma, Norman OK 73019, USA 
\email{\{khaled.abdelaal,richard.m.veras\}@ou.edu}}
\maketitle              
\begin{abstract}
    The performance of sparse matrix computation highly depends on the matching of the matrix format with the underlying structure of the data being computed on. Different sparse matrix formats are suitable for different structures of data. Therefore, the first challenge is identifying the matrix structure before the computation to match it with an appropriate data format. The second challenge is to avoid reading the entire dataset before classifying it. This can be done by identifying the matrix structure through samples and their features. Yet, it is possible that global features cannot be determined from a sampling set and must instead be inferred from local features. To address these challenges, we develop a framework that generates sparse matrix structure classifiers using graph convolutional networks. The framework can also be extended to other matrix structures using user-provided generators. The approach achieves 97\% classification accuracy on a set of representative sparse matrix shapes.

\end{abstract}

\section{Introduction}
\vspace{-1em}
\begin{figure}[ht]
    \begin{subfigure}[b]{.25\textwidth}
        \centering
        \includegraphics[width=\linewidth]{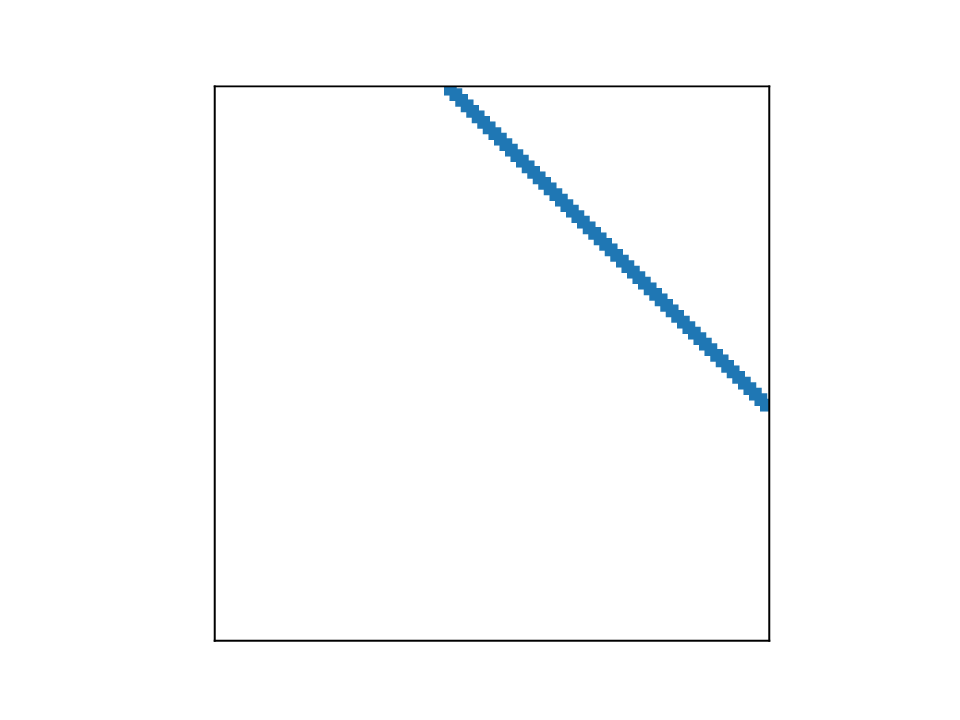}
        \caption{}
        \label{fig:diag-orig}
    \end{subfigure}%
    \begin{subfigure}[b]{.25\textwidth}
        \centering
        \includegraphics[width=\linewidth]{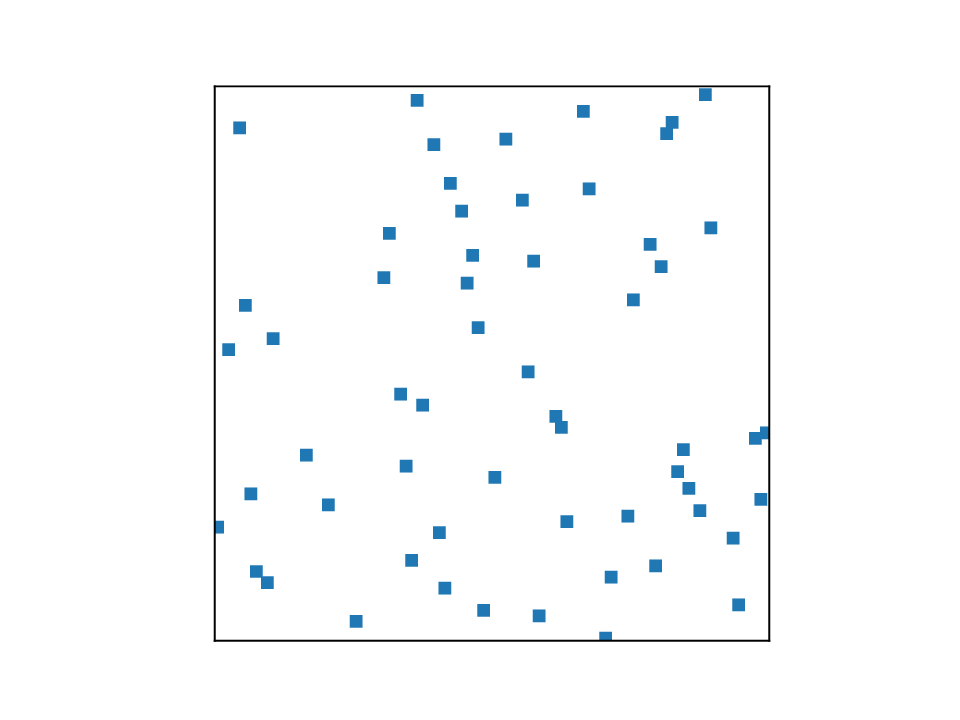}
        \caption{}
        \label{fig:diag-scrm}
    \end{subfigure}%
    \begin{subfigure}[b]{.5\textwidth}
        \centering
        \resizebox{\linewidth}{!}{%
                \begin{tabular}{c|c|c|c|c}
                     & \tiny{Diagonal} & \tiny{Random} & \tiny{Rand+Diag} & \tiny{Kronecker} \\
                    \hline
                    \tiny{Diagonal} & \textbf{500} & 0 & 0 & 0 \\
                    \hline
                    \tiny{Random} & 0 & \textbf{490} & 10 & 0 \\
                    \hline
                    \tiny{Rand+Diag} & 0 & 40 & \textbf{460} & 0 \\
                    \hline
                    \tiny{Kronecker} & 0 & 20 & 0 & \textbf{480} \\
                    \hline
                \end{tabular}
        }
        \caption{}
        \label{fig:conf-mat}
    \end{subfigure}%

    \caption{Efficacy of our classifier framework at determining structure when the data is permuted. (a) is an off-diagonal matrix, (b) is a re-labelled variant of (a), and (c) is the confusion matrix for the classifier framework on re-labelled matrices similar to (b). By using GCNs our approach is invariant to node labelling and achieves around 97\% accuracy.}
    \label{fig:de}
\end{figure}
\vspace{-1.5em}
Sparse matrices represent a fundamental building block used throughout the field of scientific computing in applications, such as graph analytics, machine learning, fluid mechanics, and finite element analysis \cite{snapnets,10.1145/2049662.2049663}.  
Such matrices appear as operands in numerous fundamental computational kernels such as sparse matrix-vector multiplication (SpMV), Cholesky factorization, LU factorization, sparse matrix -dense matrix multiplication, and matricized tensor times Khatri-Rao product (MTTKRP) among others. Building efficient algorithms for this class of kernels mainly depends on the storage format used for the sparse matrix as observed in different studies \cite{10.1145/1693453.1693471,6375570}. A variety of such formats are proposed in literature \cite{10.1145/3017994,7036061}. Hence, it is crucial to identify the structure of the matrix to choose the ideal sparse format, and eventually tailor the algorithm to that format to optimize the workload performance. 
However, identifying the structure of the matrix is not always trivial. Figure \ref{fig:diag-orig} shows a spy plot of an off-diagonal sparse matrix, and Figure \ref{fig:diag-scrm} shows the same matrix, with some of the original row indices and column indices re-labelled. It is less obvious for the latter figure to provide an insight of the original structure of the non-zeros within the matrix. Additionally, in case of huge sparse matrices, we might only have access to samples of the matrix. This could be because of computational or storage restrictions, or missing data. In these two cases (re-labelling and sub-sampling), we need efficient techniques to recognize the shape of the input matrix.

To tackle mentioned issues, we propose a framework to identify sparse matrices structures, using graph neural networks. Figure \ref{fig:conf-mat} shows the confusion matrix for the proposed framework using four sample classes on re-labelled variants. The framework design is modular, allowing users to easily augment it with new structures generators or feature sets. The main contributions of this paper are as follows:
\begin{itemize}
    \item Proposing a novel, modular Graph Neural Network framework to accurately predict the shapes of sparse matrices, including partial samples, and re-labelled variants of original matrices.
    \item Presenting a new balanced synthetic dataset for structured sparse matrices.
    \item Providing a performance analysis of graph-level classification on sparse matrices, using different feature sets. 
    \item Introducing two new compact and efficient feature sets for matrices as graphs, namely: Linear and Exponential Binned One-Hot Degree Encoding.
\end{itemize}
The rest of this paper is organized as follows: Section \ref{sec:bg} introduces the necessary background, Section \ref{sec:methods} details the design of the proposed framework, Section \ref{sec:analysis} discusses the evaluation and results of the framework, while Section \ref{sec:related} describes related work. Finally, Section \ref{sec:conclusion} summarizes the findings of the paper.

\section{Background}
\label{sec:bg}
\subsection{Graph Neural Networks}
Graph neural networks (GNNs) \cite{4700287} are a class of deep learning models that operate on graphs or networks.
Unlike traditional neural networks that operate on structured data such as images or sequences, GNNs can handle arbitrary graph structures with varying node and edge attributes, enabling them to learn powerful representations of graph-structured data.
The key idea behind GNNs is to iteratively update node embeddings by aggregating information from the embeddings of their neighbors through the "graph convolution" operation. By stacking multiple layers of graph convolution and non-linear activation functions, GNNs can learn hierarchical representations of the graph that capture both local and global information.

\subsection{Structured Matrices}
Several common structures are observed in sparse matrices, such as:

\noindent\textbf{Diagonal} all non-zeros are located on the main or a secondary diagonal. This structure represents a 1D mesh and commonly appears in various scientific and engineering applications.

\noindent\textbf{Random} the non-zero elements are randomly distributed across the matrix, with variable density. Such matrices have no specific identifiable structure. 

\noindent\textbf{Kronecker Graphs} \cite{JMLR:v11:leskovec10a} are a class of synthetic graphs that have been widely used to model real-world networks, and are generated by recursively applying the Kronecker product of a small base graph with itself. Let $A$ and $B$ be two matrices. Then, their Kronecker product $A \otimes B$ is given by
\begin{equation}
A \otimes B = 
\begin{pmatrix}
a_{11}B & \cdots & a_{1n}B \\
\vdots & \ddots & \vdots \\
a_{m1}B & \cdots & a_{mn}B
\end{pmatrix}
\end{equation}
\noindent where $a_{ij}$ are the entries of $A$. 
We use these three classes of structures, and a combination of them, as a representative set that can be combined to form more complex relationships \cite{1386651,van1992computational}. Our framework is not limited to only these structures, and they serve as an example to evaluate its performance.
\begin{figure}[ht]
    \centering
    \begin{subfigure}[b]{.5\textwidth}
        \centering
        \includegraphics[width=.8\textwidth]{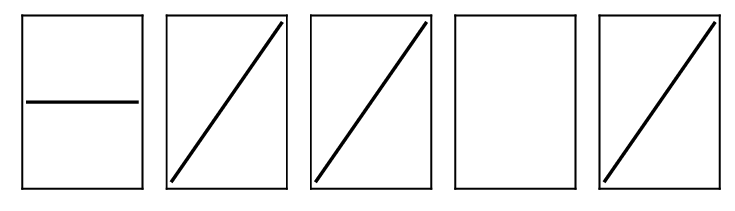}
        \caption{Diagonal}
        \label{fig:diag-deg}
    \end{subfigure}%
    \hfill
    \begin{subfigure}[b]{.5\textwidth}
        \centering
        \includegraphics[width=.8\textwidth]{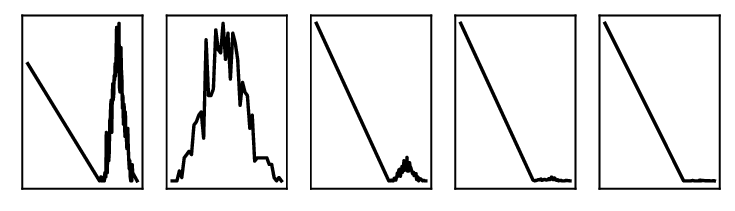}
        \caption{Random}
        \label{fig:rand-deg}
    \end{subfigure}%
    
    \begin{subfigure}[b]{.5\textwidth}
        \centering
        \includegraphics[width=.8\textwidth]{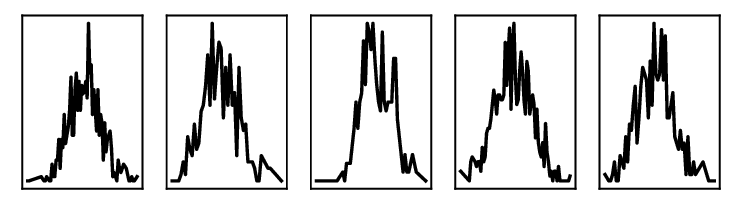}
        \caption{Random+Diagonal}
        \label{fig:rand-diag-deg}
    \end{subfigure}%
    \hfill
    \begin{subfigure}[b]{.5\textwidth}
        \centering
        \includegraphics[width=.8\textwidth]{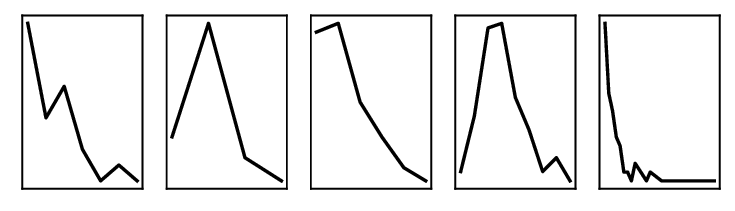}
        \caption{Kronecker Graph}
        \label{fig:rand-deg}
    \end{subfigure}%
    \caption{Global Degree Distribution for samples in each matrix (graph) class studied in this paper. In our approach we classify the shape based on local  views from sampled data.}
    \label{fig:deg-dist}
\end{figure}

\noindent\textbf{Degree as a representative node feature} Figure \ref{fig:deg-dist} illustrates that one can accurately distinguish between the different classes based on the degree distribution of the representative graph. For example, for Diagonal matrices (Figure \ref{fig:diag-deg}), the degree for all nodes is low, and is either constant or linear across all nodes. Kronecker graphs follow a power-law degree distribution, with only a few nodes having many connections (high degree) and most of the nodes having relatively few connections (low degree).
However, only the local per-node degree view may be immediately available, and not the global graph view. An example of such a case is only having a sample of the graph and not the entire graph due to storage or computational limitations. The power of GNNs can be leveraged to carry out the required task: the prediction of the sample matrix structure.  
\section{Framework Description}
\label{sec:methods}
\begin{figure*}
    \centering
    \includegraphics[width=\textwidth]{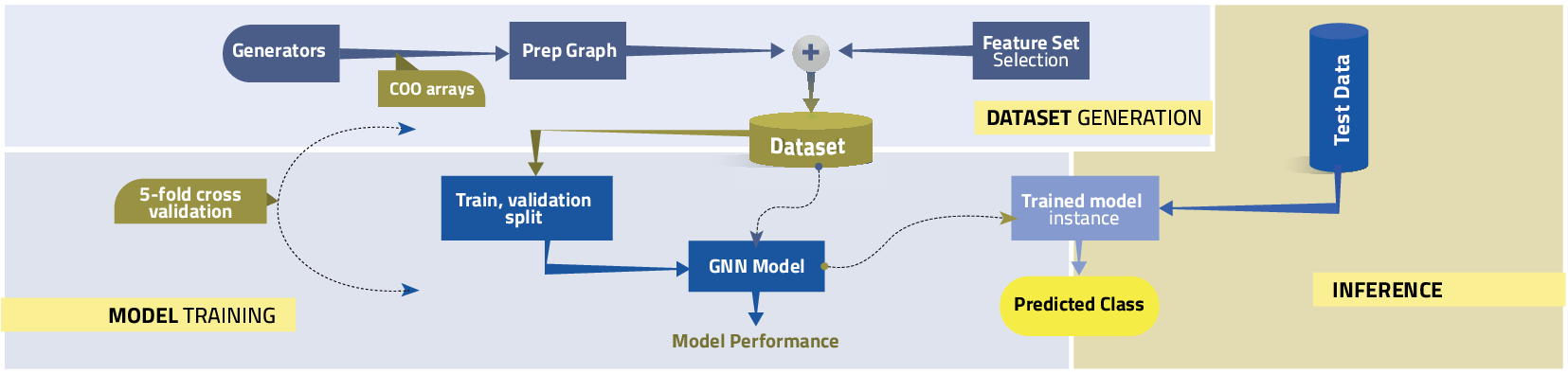}
    \caption{High-Level overview of the framework. It consists of three main phases: \textbf{dataset generation}, where the synthetic sparse matrices are generated, prepared as graphs, and have feature set attached. Then, the GNN \textbf{model training} using 5-fold cross validation to capture the model performance, and then generate a trained model instances, that is used later in the \textbf{inference} phase.} 
    \label{fig:high-level}
\end{figure*}
\noindent The goal of the proposed framework is to predict the structure of the input sparse matrix through its classification into one of the configured target classes. We use diagonal, random, diagonal+random, and Kronecker graph as examples of these classes to evaluate the performance of the framework. New structures can be seamlessly integrated. Figure \ref{fig:high-level} shows a high-level overview of the proposed framework. It consists of three main stages: Dataset generation, Model Training, and Inference. A synthetic dataset is generated using different generators for different shapes of matrices, which are then represented as graphs. In the training phase, we use GNN with 5-fold cross validation to evaluate the model performance. Finally, the trained model instance is used for later inference. 
\subsection{Dataset Generation}
\label{sec:gen}
We generate a balanced dataset of 40K synthetic sparse matrices, covering the four sample classes through individual generators. Each of the generators returns a Coordinate (COO) representation for the matrix, excluding the actual non-zero values. The COO representation is then used as the adjacency list to build the graph representation.
\subsection{Feature Set Selection}
\label{sec:feat}
A per-node feature vector is necessary for the graph neural network to classify matrices. Node degree can be calculated for rows/columns in input matrices from their graph representation.

\noindent\textbf{One-Hot Degree Encoding}
uses a number of features equal to the maximum degree + 1. A limitation of this encoding is that it requires the knowledge of the maximum degree in the entire dataset before training. Also, the required storage is proportional to the maximum degree recorded in the dataset, which increases memory requirements and reduces maximum possible batch size during training. Moreover, it poses complications during inference if the input matrix has a degree greater than the maximum degree in the training set. In our dataset, the maximum training set degree is 7710, so the length of one-hot encoding feature vector per node is 7711, limiting the maximum batch size on GPU to only 1 graph.

\noindent\textbf{Local Degree Profile}
(LDP) \cite{cai2018simple} captures the local structural information of nodes in their immediate neighborhood. LDP is calculated for each node as a five-feature vector: the node degree, the minimum degree of its neighbors, the maximum degree of its neighbors, the mean degree for its neighbors, and the standard deviation of the degrees of each neighbors. LDP features are easy to compute for any given graph. Additionally, the number of features per node is fixed, regardless of the used training data. LDP incurs low storage overhead.

\noindent\textbf{Linear Binned One-Hot Degree Encoding}
\label{sec:linear-onehot}
(LBOH) We implement a modified version of one-hot encoding, to address its limitations. LBOH works by having a fixed number of buckets for representing one-hot degrees. 
Buckets ranges are designed as follows: a set individual sequential buckets from 0 t $\alpha$ (inclusive) where $\alpha$ is a small integer (less than 10). Then, we add a set of buckets with fixed step $\beta$ from $\alpha$: $(\alpha + \beta)$, $(\alpha + 2\beta)$, ..., $(\alpha + k\beta)$ where $(\alpha + k\beta)$ is the maximum degree threshold. Any degree greater than $(\alpha + k\beta)$ is mapped to the final bucket. 
\begin{figure}
    \centering
    \includegraphics[scale=0.35]{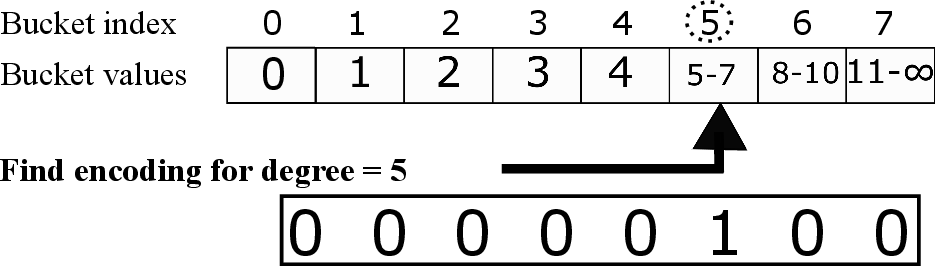}
    \caption{An example of finding the linear binned one-hot degree encoding for a node with degree = 5, where the parameters for the encoding scheme are $\alpha=5$, $\beta=3$, and $k=2$. Degree 5 is mapped to its associated bucket (5 to 7), then the bucket index (5) is represented using one-hot encoding (1 at the position where the value 5 exists, 0 otherwise).}
    \label{fig:linear-onehot}
\end{figure}
\vspace{-1.5em}
Figure \ref{fig:linear-onehot} shows an example of LBOH encoding. As opposed to One-Hot Encoding, LBOH provides a fixed number of features regardless of the maximum degree in the training dataset. At the inference stage, only the values of $\alpha$, $\beta$, and $k$ are needed. 

\noindent\textbf{Exponential Binned One-Hot Degree Encoding (EBOH)}
The main difference between EBOH and LBOH is the kind of step between buckets ranges. Instead of a linear step in LBOH, EBOH uses an exponential step to cover more degree values with a small number of features. First, the value of $\alpha$ is chosen such that  $1 \leq \alpha \leq 3$. Then for the buckets, a sequential one-to-one mapping is performed for values 0 through $2^\alpha$. For the following buckets, the upper bound (inclusive) is $2^{\alpha + i}$ where $i \in [1,k]$ and $k \in \mathbb{N}^* $. 
\begin{figure}
    \centering
    \includegraphics[scale=0.35]{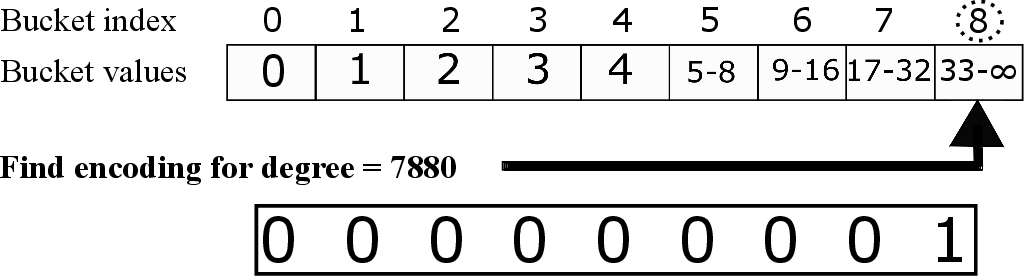}
    \caption{An example of finding the exponential binned one-hot degree encoding for a node with degree = 7880, where the parameters for the encoding scheme are $\alpha=2$ and $k=3$. Degree 7880 is mapped to its associated bucket (33 to $\infty$), then the bucket index (8) is represented using one-hot encoding (1 at the position where the value 8 exists, 0 otherwise).}
    \label{fig:expo-onehot}
\end{figure}
\vspace{-2em}
Figure \ref{fig:expo-onehot} shows an example of EBOH encoding. EBOH encoding still provides the benefit of having the number of features independent of the maximum degree in the training set. 
\subsection{The Graph Neural Network Architecture}
\vspace{-1em}
\begin{figure}[h]
    \centering
    \resizebox{0.5\linewidth}{!}{
    \begin{tikzpicture}[node distance=0.5cm, every node/.style={font=\sffamily}, scale=0.9]
    
        \node[minimum width=1cm, minimum height=0.3cm, rotate=90] (input) {\tiny Input};
        
        \node[draw, thick, rounded corners=2pt, minimum width=0.75cm, minimum height=0.3cm, rotate=90, right=of input.south, anchor=center] (conv1) {\tiny Graph Convolution};
        \node[draw, thick, minimum width=0.75cm, minimum height=0.3cm, rotate=90, right=of conv1.south, anchor=center] (relu1) {\tiny ReLU};
        \node[draw, thick, rounded corners=2pt, minimum width=0.75cm, minimum height=0.3cm, rotate=90, right=of relu1.south, anchor=center] (conv2) {\tiny Graph Convolution};
        \node[draw, thick,  minimum width=0.75cm, minimum height=0.3cm, rotate=90, right=of conv2.south, anchor=center] (relu2) {\tiny ReLU};
        \node[draw, thick, rounded corners=2pt, minimum width=0.75cm, minimum height=0.3cm, rotate=90, right=of relu2.south, anchor=center] (conv3) {\tiny Graph Convolution};
        
        \node[draw, thick, minimum width=1cm, minimum height=0.3cm, rotate=90, right=of conv3.south, anchor=center] (pooling) {\tiny Global Mean Pooling};
        \node[draw, thick, minimum width=0.75cm, minimum height=0.3cm, rotate=90, right=of pooling.south, anchor=center] (dropout) {\tiny Dropout};

        \node[draw, thick, rounded corners=2pt, minimum width=1cm, minimum height=0.3cm, rotate=90, right=of dropout.south, anchor=center] (linear) {\tiny Linear};

        \node[minimum width=1cm, minimum height=0.3cm, rotate=90, right=of linear.south, anchor=center] (output) {\tiny Output};
        
        \draw[-latex, thick] (input) -- (conv1);
        \draw[-latex, thick] (conv1) -- (relu1);
        \draw[-latex, thick] (relu1) -- (conv2);
        \draw[-latex, thick] (conv2) -- (relu2);
        \draw[-latex, thick] (relu2) -- (conv3);
        \draw[-latex, thick] (conv3) -- (pooling);
        \draw[-latex, thick] (pooling) -- (dropout);
        \draw[-latex, thick] (dropout) -- (linear);
        \draw[-latex, thick] (linear) -- (output);

    \end{tikzpicture}
    }
    \caption{Graph Neural Network architecture.}
    \label{fig:GCN}
\end{figure}
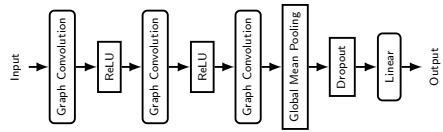
\vspace{-2em}
\noindent To identify the structure of the input matrix, the matrix is viewed as the adjacency list of a graph, enabling the use of machine learning methods designed for graph data. GNNs provide additional benefits such as allowing the use of matrices (graphs) of arbitrary sizes as input,
Also, GNNs are agnostic to node ordering. This powerful property enables re-labelling or permuting nodes in a graph representing a sparse matrix, while maintaining accurate predictions.
The machine learning task of interest is graph-level prediction since a single label (class) is needed for an entire graph (matrix). The GNN architecture is shown in Figure \ref{fig:GCN}. The hidden layers are three graph convolution layers and one linear (output) layer.  
Graph convolution is an operation where node embeddings are iteratively generated as the aggregations from the node neighborhoods. This operation is used to capture complex features of the graph. The first convolution layer aggregates information from the local neighborhood of each node. This operation is repeated in subsequent convolution layers in order to propagate information to increasingly larger neighborhoods. By the end of three convolution layers, the model has learned a hierarchical representation of the graph, where the features at each layer capture increasingly complex structural patterns. The learned representation so far is "node embeddings". Then, learned node embeddings are reduced into a single graph embedding using a global mean pooling operation (called readout layer). Samples are randomly dropped out to reduce overfitting. Finally, a linear classifier is applied to the graph embedding.
\vspace{-1em}
\begin{figure}[ht]
\centering
\begin{lstlisting}[language=Python]
def generateDiagRandom(size, threshold=2):

    """ A function to generate a Diag+Random square matrix """
    
    tuples = [(x,y) for x in range(size) for y in range(size) if (random.randint(0,10) <= threshold or x == y)]

    # seperate tuples into two lists: the row array and the column array
    coo_rep = list(map(list, zip(*tuples)))

    return coo_rep[0], coo_rep[1], [size, size]  
\end{lstlisting}
\begin{lstlisting}[language=Python]
def process(self):
    catMap = [...., {
            # Number of instances to generate for this class 
            'num_iter':10000,
            
            # Name of the generator function
            'generator':generateDiagRandom,
            
            # A string list of required generator function param 
            'gen_params':['random.randint(MIN_DIM_SIZE,MAX_DIM_SIZE)'] }]  
\end{lstlisting}
\vspace*{-1.2em}
\caption{
  The two steps needed to add a new class to the classifier framework. First (top), create a new generator function in the generators file, and second (bottom), add a dictionary entry to \lstinline{catMap} list in the \lstinline{process} method of the dataset class.
}
\label{fig:addClass}
 \vspace{-2em}
\end{figure}

\noindent\textbf{Modularity}
New shapes of matrices can be easily integrated in our framework. To achieve this, two steps are needed as shown in Figure \ref{fig:addClass}: (1) write a generator for that new shape, and (2) add an entry to the categories (shapes) map in the dataset class for this shape, containing the number of dataset instances to generate, the name of the generator function, and the different required parameters. The generator is required to return the COO representation excluding values, and the matrix dimensions. 
After generating the new dataset instances for this class, one does not need to re-train the entire model again. Transfer learning \cite{weiss2016survey} can be used to replace the last layer of the trained model with a new layer that has the appropriate number of outputs, after introducing the new shape(s). Then, the weights of all previous layers are frozen and only the new layer is trained. 
Another aspect of modularity in our framework is the ability to seamlessly attach different feature sets. Feature sets are only computed when the graph is queried. To implement a new feature set, a modification to the \lstinline{get} method of the dataset is needed. This method first reads in the graph file from disk, calculates the new feature set, and attaches it to the graph. 
\section{Analysis}
\label{sec:analysis}
We run a set of experiments to evaluate the accuracy of our approach in detecting structures, using the different feature representation discussed in Section \ref{sec:feat}. 
\subsection{Evaluation}
\noindent\textbf{Experimental Setup}
Table \ref{tab:experimental-setup} shows the experimental setup and learning parameters used in the experiments. We use PyTorch Geometric \cite{Fey/Lenssen/2019} for the GNN.
\vspace{-2.7em}
\begin{table}[ht]
    \caption{Experimental setup and Training parameters used in the experiments. \tiny{* Batch size used for traditional one-hot encoding is 1.}}
    \centering
    \begin{tabular}{lc||lc}
        \toprule
        \textbf{Component} & \textbf{Specification} & \textbf{Parameter} & \textbf{Value}\\
        \cmidrule(lr){1-2} \cmidrule{3-4}
        GPU & NVIDIA RTX A6000 & Optimizer & Adam  \\ 
        GPU Memory & 48 GB GDDR6 & Learning Rate & 0.01 \\
        CUDA Version & 11.8 & Error Criterion & Cross Entropy\\
        Main Memory & 64 GB DDR4 & Batch Size & 256{*} \\
        Operating System & Ubuntu 22.04 & Cross Validation Folds & 5 \\
        \bottomrule
    \end{tabular}
    \label{tab:experimental-setup}
\end{table}
\vspace{-1.5em}
\noindent\textbf{Evaluation Metrics}
\label{sec:eval-metrics}
To evaluate the prediction accuracy of the framework, four derived metrics are used: accuracy, precision, recall, and F1-score. We report per-class and overall accuracy and F1 score numbers, since the latter is the harmonic mean of precision and recall. Using both accuracy and F1-score helps provide a more comprehensive evaluation of the framework's performance. Accuracy gives an overall view of how well the classifier is performing, while the F1-score provides insights into its ability to correctly classify positive instances. 
\subsection{Results} %
\label{sec:results}
\vspace{-2.1em}
\begin{table*}[htbp]
    \centering
    \caption{Performance of the classifier for different degree representations}
    \label{tab:performance}
    \resizebox{\linewidth}{!}{
    \begin{tabular}{lcccccccc}
        \toprule
        & \multicolumn{8}{c}{Degree Representation} \\
        \cmidrule(lr){2-9}
        & \multicolumn{2}{c}{One-Hot Encoding} & \multicolumn{2}{c}{LDP} & \multicolumn{2}{c}{LBOH} & \multicolumn{2}{c}{EBOH} \\
        \cmidrule(lr){2-3} \cmidrule(lr){4-5} \cmidrule(lr){6-7} \cmidrule(lr){8-9}
        & Accuracy & F1 Score & Accuracy & F1 Score & Accuracy & F1 Score & Accuracy & F1 Score \\
        \midrule
        Diagonal & 1.0 & 0.90 & 1.0 & 0.97 & 1.0 & 1.0 & 1.0 & 1.0 \\
        Random & 0.90 & 0.91 & 0.64 & 0.76 & 0.92 & 0.95 & 0.95 & 0.96 \\
        Random+Diagonal & 0.86 & 0.99 & 0.98 & 0.83 & 0.96 & 0.94 & 0.97 & 0.96 \\
        Kronecker & 0.90 & 0.94 & 0.90 & 0.93 & 0.98 & 0.97 & 0.96 & 0.98 \\
        \midrule
        Overall & 0.90 & 0.90 & 0.88 & 0.88 & 0.97 & 0.97 & 0.97 & 0.98 \\
        \bottomrule
    \end{tabular}
    }
\end{table*}
\vspace{-1em}
\noindent\textbf{Classification Performance}
Table \ref{tab:performance} shows the accuracy and F1 score for the classifier using the different degree representations discussed in Section \ref{sec:feat}. Performance results show that both LBOH, and EBOH provide high prediction accuracy of around 97\% and a F1 score of around 98\%. On the other hand, traditional one-hot encoding exhibits a lower accuracy of around 90\%. One-Hot Encoding requires a significantly large number of features per node (7711), limiting the training batch size on the A6000 GPU to only one graph. This forces the optimizer to adjust the neural network weights very frequently, hence hurting the overall accuracy. 
Using LDP as a feature set exhibits variant model performance across folds depending on the validation set being used. In some folds, LDP provides high accuracy of around 97\% to 98\% similar to EBOH. In other folds, LDP fails to converge to an acceptable loss value, and ends up with an accuracy of around 74\% on the last few epochs. This performance variance across folds deems LDP unfit for the purposes of our application. It significantly fails in two classes: Random and Kronecker. It predicts Random matrices as Random+Diagonal for more than 32.5\% of the instances. This is likely due to the prevalence of the local degree neighbor summary features (the last four LDP features) instead of focusing on the node degree. This eventually results in failing to discover the global hierarchical structures in the matrix. LDP still shows perfect accuracy in case of diagonal matrices since almost all nodes in the matrix's graph have the same degree. LDP prediction quality for Kronecker graphs is also lower than other evaluated feature sets (around 81\% in some folds) for the same reasons.
 \vspace{-1.5em}
\begin{figure}[ht]
    \centering
    \begin{subfigure}{.49\textwidth}
        \centering
        \includegraphics[width=\textwidth]{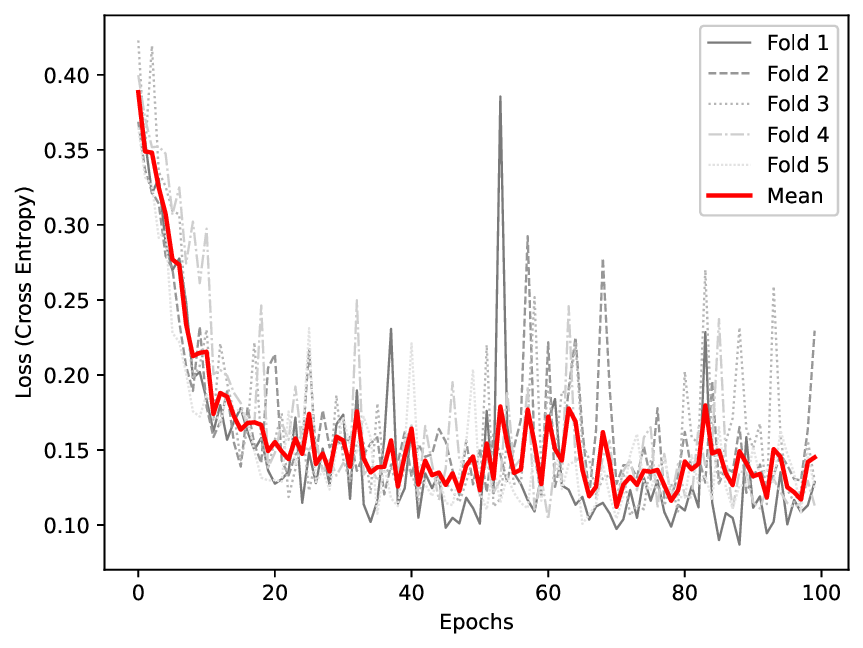}
        \caption{}
        \label{fig:loss-folds}
    \end{subfigure}
    \hfill
    \begin{subfigure}{0.49\textwidth}
        \centering
        \includegraphics[width=\textwidth]{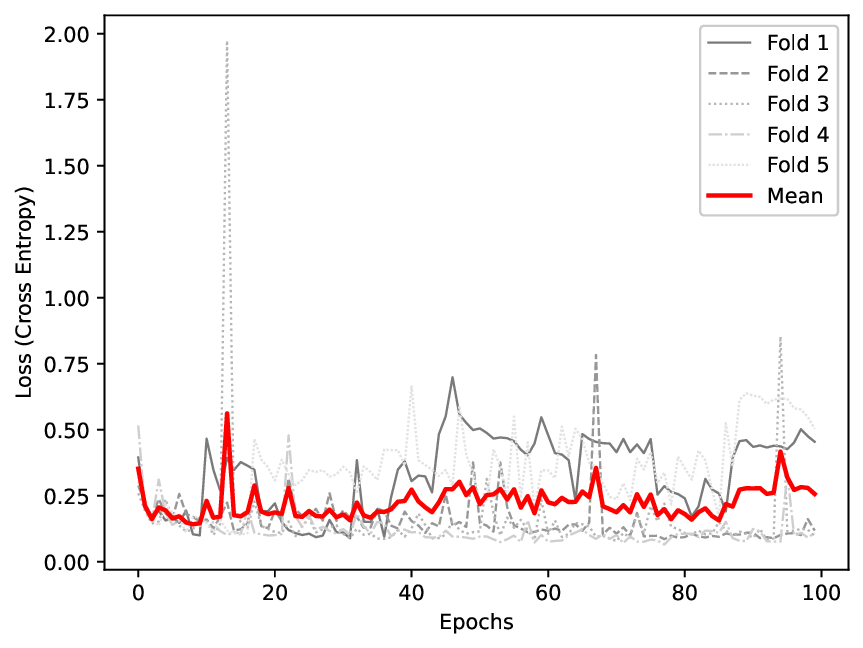}
        \caption{}
        \label{fig:loss-ldp}
    \end{subfigure}
    \caption{Cross Entropy Loss across different folds in 5-fold cross validation training using (a) EBOH, and (b) LDP feature set.}  
\end{figure}

\vspace{-2.em}

Figure \ref{fig:loss-folds} demonstrates the validation loss across the 5 different folds for EBOH. It shows almost no variance in the loss across the different folds, indicating the stability of the model's performance across folds. On the other hand, Figure \ref{fig:loss-ldp} shows the validation loss for LDP and illustrates that the loss does not converge in 2 out of 5 folds.

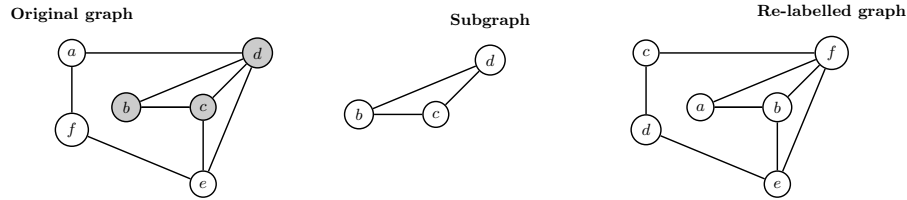
\begin{figure}[ht]
  \centering
  \resizebox{\linewidth}{!}{
  \begin{tikzpicture}[node distance=1.5cm, thick, scale=0.8]
        \node[circle,draw] (n1) at (0, 1.5){$a$};
        \node[circle,draw,fill=gray!40, below right of=n1] (n2) {$b$};
        \node[circle,draw,fill=gray!40, right of=n2] (n3) {$c$};
        \node[circle,draw,fill=gray!40, above right of=n3] (n4) {$d$};
        \node[circle,draw,below of=n3] (n5) {$e$};
        \node[circle,draw,below of=n1] (n6) {$f$};

        \draw (n2) -- (n3) -- (n4) -- (n1);
        \draw (n1) -- (n6) -- (n5);
        \draw (n2) -- (n3) -- (n5) -- (n4) -- (n2);
        \node[above=0.2cm of n1] {\textbf{Original graph}};

    \begin{scope}[xshift=7cm]
      \node[circle,draw] (n2) {$b$};
      \node[circle,draw,right of=n2] (n3) {$c$};
      \node[circle,draw,above right of=n3] (n4) {$d$};
      \draw (n2) -- (n3) --  (n4) -- (n2);
      \node[above=0.2cm of n4] {\textbf{Subgraph}};
    \end{scope}

    \begin{scope}[xshift=14cm]
        \node[circle,draw] (n1) at (0, 1.5){$c$};
        \node[circle,draw, below right of=n1] (n2) {$a$};
        \node[circle,draw, right of=n2] (n3) {$b$};
        \node[circle,draw, above right of=n3] (n4) {$f$};
        \node[circle,draw,below of=n3] (n5) {$e$};
        \node[circle,draw,below of=n1] (n6) {$d$};
        \draw (n2) -- (n3) -- (n4) -- (n1);
        \draw (n1) -- (n6) -- (n5);
        \draw (n2) -- (n3) -- (n5) -- (n4) -- (n2);
      \node[above=0.2cm of n4] {\textbf{Re-labelled graph}};
      \end{scope}
  \end{tikzpicture}
  }
  \caption{Example of generating a random sub-sample and a re-labelled variant of a graph. The original graph (left) contains six nodes. Using URNS, a random subgraph (middle) of three nodes is generated. A re-labelled variant (right) is generated using a random 1:1 mapping between the original and new node labels.}
  \label{fig:subgraph-example}
\end{figure}

\noindent\textbf{Classifying Sub-samples and Re-labelled Subgraphs}
To test the efficacy of GNN on both aspects, we generate 200 new matrices: 50 for each of the four classes. For each of them, we generate 10 subgraphs and 10 re-labelled variants. To generate the subgraphs, we use uniform random node sampling (URNS) \cite{10.1145/1150402.1150479}: nodes are randomly selected with uniform probability, as well as the edges connecting the selected nodes. Re-labelling of a graph $G$ simply renames the nodes $V$ of the graph, and produces a new graph $G'$ with the same size and degree distribution of the original graph $G$. Figure \ref{fig:subgraph-example} shows an example of both URNS and random re-labelling.

\vspace{-2em}
\begin{table}[h]
\caption{Accuracy comparison for node sampling, node re-labelling, and original graphs using EBOH feature set.}
\centering
\begin{tabularx}{\linewidth}{|c|>{\centering\arraybackslash}X|>{\centering\arraybackslash}X|>{\centering\arraybackslash}X|}
\hline
\textbf{Class} & \textbf{Node Sampling} &  \textbf{Node Re-labelling}  & \textbf{Original Graphs}\\ \hline
Diagonal & 1  & 1 & 1 \\ \hline
Random & 0.83  & 0.98 & 0.98 \\ \hline
Random+Diagonal & 0.92  & 0.92 & 0.92 \\ \hline
Kronecker & 0.94  & 0.96 & 0.96\\ \hhline{|====|}
\textbf{Overall} & \textbf{0.92}  & \textbf{0.97} & \textbf{0.97} \\ \hline
\end{tabularx}
\label{tab:matrix_accuracy}
\end{table}
\vspace{-1.5em}

 Table \ref{tab:matrix_accuracy} shows the model's performance on subgraphs and re-labelled variants as compared to original full graphs. The table shows that re-labelling node has no impact on the classification accuracy; it shows the same overall accuracy of around 97\% which is observed for the original graphs. This is expected because the arrangement of nodes in a graph is irrelevant, since the graph has the same degree distribution. For subgraphs generated using URNS of larger graphs, the overall accuracy drops to around 92\%. The reason being that random node sampling can alter the degree distribution of the graph. The random choice of nodes can result in either isolated nodes (no edges) or much lower degree nodes as compared to the original graph. This affects the accuracy specially for complex shapes such as random and Kronecker graphs. One way to reduce the accuracy loss for samples is to use a more sophisticated graph sampling technique rather than randomly selecting nodes or edges.

\vspace{-2em}
\section{Related Work}
\label{sec:related} 
\noindent\textbf{Prediction on Sparse Matrices} Several studies investigated the use of machine learning to predict the optimal sparse format for SpMV on CPU and GPU \cite{10.1145/3178487.3178495,10.1145/3218823,10.1145/2304576.2304624,6748055,10.1145/2499370.2462181,7573853}. Our framework does not directly predict the best sparse format, instead, we only predict the structure of the input matrix. This allows de-coupling the sparsity pattern from the sparse format, following the argument adopted by AlphaSparse \cite{10.5555/3571885.3571972} since our framework also allows the seamless integration of new classes. Existing techniques collect a set of features from each matrix such as: the number of diagonals, the ratio of true diagonals to total diagonals, the (maximum) number of non-zeros per row, the variation of the number of nonzeros per row, the ratio of nonzeros in DIA and ELL data structures, and a factor or power-law distribution. We only need to calculate one feature per node: its degree. Also, \cite{10.1145/3178487.3178495} uses a CNN approach to treat matrices as images, and in order to fix the size of the matrix, they normalize input matrices into fixed size blocks, losing partial matrix information in the process. In contrast, our approach handles arbitrary sizes of matrices, without losing precision, leveraging the power of Graph Neural Networks. We can optionally sample large matrices and maintain high prediction accuracy. An additional benefit to our framework is that it is order in-variant, since matrices are represented as graphs. 

\noindent\textbf{Graph Representation for Learning} Representing non-attribute graphs is an open problem \cite{10.1145/3511808.3557661}. Common approaches employ graph properties such as node degree, more specifically a one-hot encoding of the degree \cite{xu2018how}.
One-hot encoding suffers from numerous limitations (Section \ref{sec:feat}). LDP \cite{cai2018simple} provides a compact representation for graph using five features per node. Although the computation of such feature vector is efficient, using LDP results in unreliable model performance for our task (Section \ref{sec:results}).  Both our representations (LBOH and EBOH) outperform one-hot encoding and LDP while addressing their shortcomings.
\section{Summary}
\label{sec:conclusion}
In this paper, we proposed a GNN based framework to classify structured sparse matrices. We introduced two novel non-attribute graph representations based on node degrees: LBOH, and EBOH. We evaluated the efficacy of our framework on a synthetic, balanced dataset of matrices that we generated containing random matrices from four sample classes: diagonal, random, random+diagonal, and Kronecker graphs. Performance results demonstrate a high classification accuracy of 97\% for the framework when using our feature sets: LBOH and EBOH. They also show high accuracy of 92\% and 97\% on random node subsamples and re-labelled variants respectively. Our framework is modular, allowing the inclusion of additional classes with minimal user effort. Future endeavors target the automatic generation of the optimal sparse data format and algorithm for sparse matrix kernels, using the obtained prediction results from the current framework.

\bibliographystyle{splncs04}
\bibliography{refs}

\end{document}